\documentclass[twoside,11pt]{amsart} 
\usepackage{amsfonts,amssymb,latexsym,amsthm, epsfig,amsmath,oldgerm,amscd}
\newcommand{\mf}{\mathfrak}

\newcommand{\ra}{\rightarrow}
\newcommand{\Ra}{\Rightarrow}

\newcommand{\eps}{\varepsilon}

\newcommand{\mbb}{\mathbb}
\newcommand{\tn}{\textnormal}

\newtheorem{de}{Definition}[section]
\newtheorem{re}[de]{Remark}
\newtheorem{pr}[de]{Proposition} 
\newtheorem{tr}[de]{Theorem}
\newtheorem{lm}[de]{Lemma} 

\newtheorem{nt}[de]{Notation} 

\newtheorem{co}[de]{Corollary}

\newcommand{\lb}{\linebreak}

\def\vp{\rm \vspace{0.2cm}}
\def\M{\rm Max}
\def\m{\mf{m}}

\def\hb{\hfill$\Box$}
\def\aut{\rm Aut} 
\def\sp{\rm Spec} 
\def\tran{\rm Trans}
\def\fl{\rm ETrans}
\def\Um{\rm Um}
\def\GL{\rm GL}

\def\SL{\rm SL}

\def\EO{\rm EO}

\def\SO{\rm SO}
\def\E{\rm E}

\def\G{\rm G}
\def\Sp{\rm Sp}

\def\k{\rm K_1}

\def\I{\rm I}

\def\O{\rm O}
\def\ESp{\rm ESp}
\def\EO{\rm EO}

\def\Um{\rm Um}

\begin{document}
\title{The Pillars of Relative Quillen--Suslin Theory}
\author{Rabeya Basu, Reema Khanna \& Ravi A. Rao}
\date{}
\maketitle
{\small \begin{center}
\end{center}}
{\small \begin{center}{\it Key words: Bilinear forms, Symplectic and 
Orthogonal forms}
\end{center}}
{\small Abstract:  We deduce the relative version of the equivalences relating 
the relative Local Global Principle and the Normality of the relative 
Elementary subgroups of the traditional classical groups; {viz.} general 
linear, symplectic and orthogonal groups. This generalizes our previous result 
for the absolute case; {\it cf.} \cite{rrr}}.

\section{Introduction} 

The main pillars of the Horrocks--Quillen--Suslin theory were developed in the 
papers \cite{hor}, \cite{QUI}, \cite{s}. In \cite{hor}  the 
{\it Monic Inversion Principle}, in \cite{QUI} 
the {\it Local-Global Principle}, and in \cite{s} the {\it Normality of the 
Elementary subgroup ${\E}_n(R)$}, were established. In \cite{s} the ${\k}$ 
analogues 
of both the Monic Inversion Principle and the Local-Global Principle were 
developed. In addition, Suslin established the {\it Normality of the Elementary Linear 
subgroup} ${\E}_n(R)$ in the general linear group ${\GL}_n(R)$  over a module finite ring $A$, when $n \geq 3$. This 
was appeared in \cite{Tu}.

In \cite{rrr} the authors had established, for  classical linear groups, 
{\it viz.} the linear, symplectic and orthogonal groups, that the 
Quillen--Suslin's Local-Global 
Principle for the pair $({\GL}_n(R[X]), {\E}_n(R[X])$ and  Suslin's 
Normality Principle were {\it equivalent} in the sense that if one holds 
then so does the other. Recently, in \cite{rrsy} a further unification 
of these three principles was achieved.  

In this article, we develop the equivalence of a {\it relative Quillen's
Local-Global Principle} and a {\it normality of the relative elementary 
subgroup}; {\it cf.} Theorem \ref{nlg} for the precise 
equivalent statements. 

We refer the reader to the Introduction of \cite{rrr} where recent 
developments of the Quillen--Suslin theory are discussed in detail. The study 
of the relative Local-Global Principle with respect to an extended ideal 
began in \cite{hpr}; and was developed in \cite{hstep} for the Chevalley 
groups. 

The proofs of the equivalent statements in this paper are done in an 
analogous manner to that done in \cite{rrr}. This was 
possible due to a recent argument, which is detailed in \cite{jkrs}, and which 
first appeared in the thesis of Anjan Gupta \cite{ag}. This argument works with 
the Noetherian excision ring $R\oplus I$ rather than the use of the 
(non-Noetherian) Excision ring $\mbb{Z} \oplus I$, and the Excision theorem 
of W. van der Kallen in \cite{vand}, as is commonly used. We refer \cite{ag1} to see other interesting 
applications of the Noetherian Excision rings.

For the sake of being self-contained we have detailed the arguments of the 
various equivalences. However, we note that we could have alternatively 
{\it deduced} the 
implications from the corresponding implications done in \cite{rrr} via this 
Noetherian Excision ring argument. 

\section{Definitions and Notations}

Let $R$ be a commutative ring with $1$, and $I\subset R$ an ideal. We
refer \cite{rrr} for the standard definitions and facts of the
general linear, symplectic and orthogonal groups, and their elementary
subgroups. Let $\sigma$ denote the permutation of the natural
numbers given by  $\sigma(2i)=2i-1$ and $\sigma(2i-1)=2i$. With respect to this permutation we 
define following classical groups. 

For an integer $m > 0$, the symplectic  group of size $2m\times 2m$ is defined with respect to the
alternating matrix $\psi_m$ corresponding  to the standard symplectic
form $$\psi_m=\underset{i=1}{\overset{m}\sum} e_{2i-1,2i}-
\underset{i=1}{\overset{m}\sum} e_{2i,2i-1}.$$ For the orthogonal
group we have considered  symmetric matrix $\widetilde{\psi}_m$ 
corresponding to the  standard hyperbolic form 
$$\widetilde{\psi}_m=\underset{i=1}{\overset{m}\sum} e_{2i-1,2i}+ 
\underset{i=1}{\overset{m}\sum} e_{2i,2i-1}.$$ 

\begin{de} \tn{{\bf Symplectic Group ${\Sp}_{2m}(R)$:}  The group of all
non-singular $2m\times 2m$ matrices}  $ \{\alpha\in
{\GL}_{2m}(R)\,\,|\,\,\alpha^t \psi_m \alpha = \psi_m\}$.
\end{de} 

\begin{de} \tn{{\bf Orthogonal Group ${\O}_{2m}(R)$:}  The group of all
non-singular $2m\times 2m$ matrices}  $\{\alpha\in
{\GL}_{2m}(R)\,\,|\,\,\alpha^t \widetilde{\psi}_m \alpha =
\widetilde{\psi}_m\}$.
\end{de}

\begin{de} \tn{{\bf Elementary Symplectic Group ${\ESp}_{2m}(R)$:}  For
$1\le i\ne j\le 2m$ we define},
\begin{align*} se_{ij}(z) &= {\I}_{2m}+ze_{ij} \tn{  if  } i=\sigma(j) \\
&= {\I}_{2n}+ze_{ij}-(-1)^{i+j}z e_{\sigma (j) \sigma (i)} \tn{ if }
i\ne \sigma(j) \tn{ and } i<j. 
\end{align*} \tn{It is clear that  when $z\in R$ all these matrices
belong to  ${\Sp}_{2m}(R)$. We call them the elementary symplectic
matrices over $R$  and the group generated by them is called
elementary symplectic group.}
\end{de}

\begin{de} \tn{{\bf Elementary Orthogonal Group ${\ESp}_{2m}(R)$:}  For
$1\le i\ne j\le 2m$ we define},
\begin{align*} oe_{ij}(z) &=  {\I}_{2n}+ze_{ij}-z
e_{\sigma(j)\sigma(i)}\tn{ if } i\ne \sigma(j)  \tn{ and } i<j. 
\end{align*}  \tn{It is clear that  when $z\in R$ all these matrices
belong to  ${\O}_{2m}(R)$. We call them the elementary orthogonal
matrices over $R$  and the group generated by them is called
elementary orthogonal group.}
\end{de}

\begin{nt} \tn{ In the sequel ${\rm M}(n, R)$ will denote the set of all
$n\times n$  matrices,  ${\G}(n,R)$ will denote either the linear group
${\GL}_n(R)$, the symplectic group  ${\Sp}_{2m}(R)$, or the orthogonal group
${\O}_{2m}(R)$,  where $2m = n$. ${\rm S}(n,R)$ will denote either the special
linear group ${\SL}_n(R)$, the symplectic group ${\Sp}_{2m}(R)$, or the special orthogonal group ${\SO}_{2m}(R)$, when
$R$ is a  commutative ring.  
Similarly, ${\E}(n,R)$ will denote the
corresponding elementary subgroups  ${\E}_n(R)$, ${\ESp}_{2m}(R)$,
${\EO}_{2m}(R)$ respectively. To denote the generators of ${\E}(n,R)$ we
shall use the symbol $ge_{ij}(x)$, $x\in R$. }
\end{nt}
\begin{de} \tn{ The elementary subgroup ${\E}(n,I)$ with respect to the
ideal $I$ is the subgroup of ${\E}(n, R)$ generated as a group  by the
elements $ge_{ij}(x)$, for $x \in I$.  The {\it{relative elementary
group }} ${\E}(n, R,I)$ is the normal closure of ${\E}(n, I)$ in 
${\E}(n,R)$. }
\end{de} 
\begin{nt} \tn{The relative subgroups of ${\G}(n,R)$ and ${\rm S}(n,R)$ will be
denoted by ${\G}(n,R, I)$ and ${\rm S}(n,R, I)$ respectively. 
{\it i.e.} } $${\G}(n,R, I) = \{\alpha\in {\G}(n,R)\,|\,\alpha\equiv {\I}_n ~~{\rm modulo}~~ I\},$$
$${\rm S}(n,R, I) = \{\alpha\in {\rm S}(n,R)\,|\,\alpha\equiv {\I}_n ~~{\rm modulo}~~ I\}.$$
\tn{ For an ideal  $I$
in $R$, its extension in the ring $R[X]$, {\it i.e.} $I\otimes_R R[X]$
will be denoted by $I[X]$. } 

\tn{Similarly, ${\Um}_n(R,I)$ will denote the set
of all unimodular rows of length $n$ which are congruent to
$e_1=(1,0,\ldots,0)$ modulo $I$.}

\tn{We will mostly use localizations with respect to two types of
multiplicatively closed subsets of $R$. {\it viz.}
$S=\{1,s,s^2,\ldots\}$, where  $s\in R$ is a non-nilpotent, non-zero
divisor, and $S=R \setminus \mf{m}$ for $\mf{m}\in \M(R)$.  By $I_s[X]$ and
$I_{\mf{m}}[X]$ we shall mean the extension of $I[X]$ in $R_s[X]$ and
$R_{\mf{m}}[X]$ respectively.}

\end{nt}  {\bf Blanket Assumption:}  We assume that $n\ge 3$, when
dealing with the linear case and $n = 2m$, with $m\ge 2$, when considering
the symplectic and orthogonal  cases. While dealing with the
orthogonal groups we shall consider only isotropic vectors; {\it i.e.}
all such non-zero vectors which are orthogonal to themselves  with respect
to the given non-degenerate bilinear form. Throughout the article we
shall assume  2 is invertible in the ring $R$.
\begin{nt}  \tn{For any column vector $v\in R^n$ we denote by
$\widetilde{v}=v^t.\psi_n$ in  the symplectic case and
$\widetilde{v}=v^t.\widetilde{\psi}_n$ in the orthogonal case.}
\end{nt} 
\begin{de} \tn{We define the map $M:R^n \times R^n \ra M(n,R)$ and the
inner product $\langle ~,\rangle$ as follows: Let $v, w$ be column vectors in $R^n$. Then,}
\begin{align*} M(v,w) & = v.w^t, \,\,\,\,  \tn{ when dealing with the
case }  {\G}(n,R)={\GL}_n(R). \\ & = v.\widetilde{w}+w.\widetilde{v},
\,\,\,\, \tn{ when }  {\G}(n,R)={\Sp}_{2m}(R).\\ & =
v.\widetilde{w}-w.\widetilde{v}, \,\,\,\,  \tn{ when }
{\G}(n,R)={\O}_{2m}(R).\\  \langle v,w\rangle & = v^t.w, \,\,\,\, \tn{ when
} {\G}(n,R)={\GL}_n(R). \\ & = \widetilde{v}.w, \,\,\,\, \tn{ when }
{\G}(n,R)={\Sp}_{2m}(R) \tn{ or } O_{2m}(R).
\end{align*} 
\end{de} 
\begin{nt} \tn{For any $\alpha\in {\G}(n,R)$, as usual $\alpha \perp {\I}_r$
denotes its  embedding in ${\G}(n+r,R)$, where $r$ is even for non-linear
cases.}
\end{nt}

To deduce the relative case from the absolute case we consider the
``Noetherian Excision ring''. 
\begin{de} {\bf (The ring $R\oplus I$)}:  \tn{ Let $I$ be an ideal in
the ring $R$. We construct the new ring $R \oplus I$ by defining
addition and co-ordinate wise  multiplication as follows:
$$(r \oplus j)(s \oplus i)=rs \oplus (sj+ri+ij) \tn{~~for~~} r, s \in R \tn{~~and~~} i, j \in I.$$
There is a natural homomorphism $\phi : R \oplus I \longrightarrow R$
given by $(r \oplus i)\rightarrow r+i \in R$. }
\end{de} 

\begin{nt}\tn{ Let ${\E}(n, I) = \{\alpha\in {\rm S}(n,R)\,|\,\alpha\equiv {\I}_n ~~{\rm modulo}~~ I\}.$ In general, ${\E}(n, I)$ is not 
normal in ${\G}(n,R)$. 
By  ${\E}(n,R,I)$ 
we mean the the normalization of ${\E}(n,I)$ in ${\G}(n,R)$,
{\it i.e.} the  relative elementary group  generated by elements of
the type  $ge_{ij}(f)ge_{ji}(h)(ge_{ij}(f))^{-1}$, where
$f \in R$ and $h \in I$. While working on the polynomial ring $R[X]$, by writing $\alpha(X) \in {\E}(n,R[X],I[X])$ we  mean $\alpha(X)$ is 
 ${\I}_n$ modulo $I$, and of the form $ge_{ij}(f(X))ge_{ji}(h(X))(ge_{ij}(f(X)))^{-1}$, where
$f(X) \in R[X]$ and $h(X) \in I[X]$, as ${\E}(n,R[X],I[X])$ 
is the normalization of ${\E}(n,I[X])$ in ${\G}(n,R[X])$.}
\end{nt}

\begin{lm} \label{elementary} If $\epsilon \in {\E}(n,R,I)$, then there
exists $\widetilde{\epsilon} \in {\E}(n,R \oplus I)$ such that $\phi
(\widetilde{\epsilon}) = \epsilon$. $($In fact, the converse is also true$)$.
\end{lm} {\bf Proof.}  Let $\epsilon = (\epsilon_{ij})$ be a generator
of the type $ge_{ij}(a)ge_{ji}(x)ge_{ij}(-a)$, where $a\in R$ and
$x\in I$. We then have $1$'s, $1+ax$ and $1-ax$ on the diagonal of
$\epsilon$ and zeros,  $-a^2x$ and $x$ as the non-diagonal elements.
We get a new matrix $\widetilde{\epsilon}$ by taking the corresponding
diagonal elements as $(1,0)$, $(1,ax)$ and $(1,-ax)$ and the
corresponding non-diagonal elements as $(0,0)$, $(0,-a^2x)$ and
$(0,x)$ which are elements of the ring $R \oplus I$.   Using the
definition of multiplication in the ring $R \oplus I$, we can see that
\begin{eqnarray*}
\widetilde{\epsilon} &=& ge_{ij}((a, 0))ge_{ji}((0, x))ge_{ij}(-(a,0))  
\in {\E}(n,R \oplus I),
\end{eqnarray*}
and  applying the
homomorphism $\phi$ to it we obtain $\phi (\widetilde{\epsilon}) =
\epsilon$. \hb

\begin{lm} Let $\alpha \in {\G}(n,R,I)$. Then there exists
$\widetilde{\alpha} \in {\G}(n,R \oplus I)$ such that $\phi
(\widetilde{\alpha}) = \alpha$.
\end{lm} {\bf Proof.}  Let $\alpha = ({\alpha}_{ij}) \in {\G}(n,R,I)$.
Then ${\alpha}_{ii}= 1 +a_{ii}$ and ${\alpha}_{ij}= a_{ij}$ for $i \neq
j$ where $a_{ij} \in I$ for all $i,j$.  We
get a new matrix $\widetilde{\alpha} = \widetilde{\alpha }_{ij}$,
where $\widetilde{\alpha }_{ii}= (u_i,a_{ii})$  and
$\widetilde{\alpha}_{ij} = (0,a_{ij})$ for $i \neq j$. The entries
in $\widetilde{\alpha} $ are in the ring  $R \oplus I$. Using the
definition of multiplication in the ring $R \oplus I$, we can see that
$\widetilde{\alpha} \in {\G}(n,R \oplus I)$  and applying the
homomorphism $\phi$ we obtain $\phi (\widetilde{\alpha}) =
\alpha$. \hb\vp 

Now we state the main theorem of this article. For the absolute case;
{\it i.e.} for $I=R$ we refer to \cite{rrr}.

\section{Equivalence: Relative L-G Principle and Normality}

\begin{tr} \label{nlg}

Let $R$ be a commutative ring with identity, and $I\subsetneq R$ an
ideal of the ring $R$.  Let $v, w$ be column vectors in $R^n$ with $w\in I^n$. Then the followings are equivalent:

\begin{enumerate}
\item [(1)] $(${\bf Normality}$)$: ${\E}(n,R,I)$ is a normal subgroup of
${\G}(n,R)$. \vp 
\item [(2)] ${\I}_{n} + M(v,w) \in E_n(R,I)$ if  $v \in {\Um}_n(R,I)$ and
$\langle v,w\rangle =0$ and $w \in I^n$. \vp 
\item [(3)] $(${\bf Local Global Principle}$)$: \\If $\alpha (X) \in
{\G}(n,R[X],I[X])$ ; $\alpha (0) = {\I}_n$ and  ${\alpha }_{\mf{m}}(X) \in
{\E}(n,R_{\mf{m}}[X], {\I}_{\mf{m}}[X])$ for all $ \mf{m} \in \M(R)$ then
$\alpha (X) \in {\E}(n,R[X],I[X])$. \vp 
\item [(4)] $(${\bf Dilation Principle}$)$: \\If $\alpha (X) \in
{\G}(n,R[X],I[X])$ ; $\alpha (0) = {\I}_n$ and  ${\alpha }_s(X) \in
{\E}(n,R_s[X], I_s[X])$ for some non-nilpotent element $s \in R,$ then
$\alpha (bX) \in {\E}(n,R[X],I[X])$ for $b \in (s^l)$, $l\gg 0$.
$($Actually,  we mean there exists some $\beta(X)\in {\E}(n,R[X],I[X])$
such that $\beta(0)={\I}_n$ and $\beta_s(X)=\alpha_s(bX)$. But, since
there is no ambiguity,  for simplicity we are using the notation
$\alpha(bX)$ instead of $\beta_s(X)).$ \vp 
\item [(5)] Let $\alpha (X) ={\I}_{n} + X^d M(v,w)$ for some integer
$d\gg 0$, $v \in {\E}(n,R,I)e_{1}$, $w\in I$ with $\langle v,w\rangle=0$. Then one
gets $\alpha (X) \in {\E}(n,R[X],I[X])$. Moreover,  $\alpha (X)$ can be
expressed as a product decomposition of the form $\Pi ge_{ij}(Xh(X))$
for $d\gg 0$ and $h(X) \in I[X]$. \vp 
\item [(6)] ${\I}_{n} + M(v,w) \in {\E}(n,R,I)$ if $v \in {\E}(n,R,I)e_{1}$, $w
\in I^n$  and $\langle v,w\rangle=0$. \vp 
\item [(7)] ${\I}_{n} + M(v,w) \in {\E}(n,R,I)$ if $v \in {\G}(n,R,I)e_{1}$, $w
\in I^n$  and $\langle v,w\rangle=0$.
\end{enumerate}
\end{tr} 

\begin{re} Since (6) will be established in Lemma \ref{key5}, it follows 
that all the above statement (1)-(7) of Theorem \ref{nlg} hold for commutative $($In fact, for almost commutative$)$ rings.  
\end{re}

\vskip0.15in
 
Before proving the theorem we first collect a few lemmas.

\begin{lm} \label{nor} The group  ${\E}(n,R, I)$ satisfies the
property:  
$$[{\E}(n,R, I), {\E}(n,R)] = {\E}(n,R, I).$$ 
\end{lm}
{\bf Proof.} {\it cf.}\cite{BA} for the general linear groups, (\cite{KOP}, Theorem 1.1) for the symplectic groups and (\cite{SUSK}, $\S 2$) for the orthogonal groups. 
\hb \vp

Below we state few useful well-known lemmas. For
the proofs  {\it cf.} \cite{BA} for the linear groups, \cite{KOP} for
the symplectic groups, \cite{SUSK} for the orthogonal groups. For a uniform proof {\it cf.} \cite{rrr}, \cite{RB1}.  
The analogue results for the relative cases follows from the proofs of the absolute cases. 
\vp 

\begin{lm} \label{key01} \tn{{\bf (Splitting Property):}}
$ge_{ij}(x+y)=ge_{ij}(x)ge_{ij}(y)$, $\forall~ x,y\in R$. 
\end{lm} 

\begin{lm} \label{key02} Let $G$ be a group, and $a_i$, $b_i \in G$,
for  $i = 1, \ldots, n$. Then   $ {\underset{i=1}{\overset{n}\Pi}}a_i
b_i= {\underset{i=1}{\overset{n}\Pi}}r_ib_ir_i^{-1}
{\underset{i=1}{\overset{n}\Pi}}a_i$, where $r_i={\underset{j=1}
{\overset{i}\Pi}} a_j$. 
\end{lm}  
\begin{lm} \label{key1} The group  ${\G}(n,R[X],(X)) \cap {\E}(n,R[X],
I[X])$ is generated by the elements of the type $ \eps
ge_{ij}(Xh(X))\eps^{-1}$, where  $\eps \in {\E}(n,R[X])$, $h(X)\in I[X]$.
\end{lm}

\begin{lm} \label{key3} For $m>0$, and $h(Y)\in I[Y]$, there are $h_t(X,Y,Z)\in I[X,Y,Z]$
such  that
$$ge_{pq}(Z)ge_{ij}(X^{2m}h(Y))ge_{pq}(-Z)= \underset{t=1}{\overset{k}\Pi}
ge_{p_tq_t}(X^m h_t(X,Y,Z)).$$ 
\end{lm} 

\begin{co} \label{key4} If $\eps=\eps_1\eps_2\cdots\eps_r$, where each
$\eps_j$ is an elementary generator, and $h(Y)\in I[Y]$, then there are $h_t(X,Y)\in
I[X,Y]$ such that
$$\eps ge_{pq}(X^{2^rm}h(Y))\eps^{-1}=\underset{t=1}{\overset{k}\Pi}
ge_{p_tq_t}(X^m h_t(X,Y)).$$
\end{co}  {\bf Proof.}  Follows by induction on $r$ and using Lemma
\ref{key3}.  \hb

We show that statement (6) of Theorem \ref{nlg} is true over an
arbitrary  associative ring $R$ with $1$.

\medskip

\begin{lm} \label{key5}  Let $R$ be a ring and $v\in {\E}(n,R,I)e_1$. Let
$w\in I^n$ be a  column vector such that $\langle v,w\rangle=0$. Then
${\I}_n+M(v,w)\in {\E}(n,R,I)$.
\end{lm} {\bf Proof.}  Let $v=\eps e_1$, where $\eps =(\eps_{ij}) \in
{\E}(n,R,I)$. Hence $\eps_{ii}=1+a_{ii}$ and $\eps_{ij}=a_{ij}$ for $i\ne
j$, where $a_{ij}\in I$ for all $i, j$.   Let
$\widetilde{\eps}=\widetilde{\eps}_{ij}$, where
$\widetilde{\eps}_{ii}=(1,a_{ii})$, and
$\widetilde{\eps}_{ij}=(0,a_{ij})$ for $i\ne j$.  Let
$\widetilde{e_1}=((1,0),(0,0),\ldots,(0,0))$, and 
$$\widetilde{v}= ((1,v_1),(0,v_2),\ldots,(0,v_n))\in (R \oplus I)^n,$$  
$$\widetilde{w}=((0,w_1),(0,w_2),\ldots,(0,w_n)) \in (0 \oplus I)^n.$$ 
Then it follows that 
$${\I}_n+M(\widetilde{v},\widetilde{w})=\widetilde{\eps} ({\I}_n+ M(\widetilde{e}_1,\widetilde{w}_1))\eps^{-1},$$  
$$\tn{~~and~~}\widetilde{w_1} = \begin{cases}
\widetilde{\eps}^t \widetilde{w} & \mbox{ for linear case } \\
\widetilde{\eps}^{-1}\widetilde{w} & \mbox{ otherwise.}
\end{cases}$$ Since $\langle (\widetilde{e_1},
\widetilde{w}_1)\rangle=\langle \widetilde{v},\widetilde{w}\rangle=0$,
we get
$$\widetilde{w_1}^t=
\begin{cases}  ((0,0),(0,w_{12}),(0,w_{13}),\dots,(0,w_{1n})) & \mbox{
for linear case }\\ ((0,w_{11}),(0,0),(0,w_{13}),\dots,(0,w_{1n})) &
\mbox{ otherwise. } 
\end{cases}$$  Therefore,
$${\I}_n+M(\widetilde{v},\widetilde{w})= \begin{cases} \underset{j=2}{\overset{n}\Pi}
\,\widetilde{\eps} \,ge_{1j}(0,w_{1j})\,\widetilde{\eps}^{-1} &
\mbox{for linear case } \\ \underset{j\ne
2}{\underset{j=1}{\overset{n}\Pi}} \,\widetilde{\eps}
\,ge_{1j}(0,w_{1j})\,\widetilde{\eps}^{-1} & \mbox{ otherwise }.
\end{cases}$$  Hence ${\I}_n+M(\widetilde{v},\widetilde{w})\in
{\E}(n,R\oplus I,0\oplus I)$.  Now applying the homomorphism $\phi$ it
follows that ${\I}_n+M(v,w)\in {\E}(n,R,I)$; as desired. \hb \vp 

Note that the above implication is true for any associative ring with
identity. 

\begin{re} \label{finite}  \tn{It is well-known that every ring is a direct 
limit of Noetherian rings. Hence we may consider $R$ to be Noetherian.}
\end{re}

We shall use following lemma frequently and sometime in a subtle way; 
{\it e.g.} for the implication $(4)\Ra (3)$.

\begin{lm} \label{noeth} {\rm ((\cite{HV}, Lemma 5.1))} Let $R$ be a Noetherian ring and $s\in R$.
Then there exists a natural number  $k$ such that the canonical homomorphism
${\G}(n,s^kR) \ra {\G}(n, R_s)$  $($induced by localization
homomorphism $R \ra R_s)$ is injective.  Moreover, it follows that
the map  ${\E}(n,R,s^kR) \ra {\E}(n,R_s)$ for $k\in \mbb{N}$ is
injective. 
\end{lm} 
{\bf Proof of Theorem \ref{nlg}:} We shall assume the
result for the absolute case; {\it i.e.} when $I=R$.  The implication
$(7)\Ra (6)$: Obvious. We prove, $(6)\Ra (5)$: \vp 

Note that we have assumed that $(6)$ holds for any commutative ring, in 
particular for the ring $R[X]$, and the matrix $I_n + XM(v, w)$. 
Replacing $R$ by $R[X]$ in $(6)$ we get that ${\I}_n+XM(v,w)\in {\E}(n,R[X],
I[X])$. Let $v=\eps e_1$, where $\eps \in {\E}(n,R,I)$.  As before, let
$\widetilde{v}= ((1,v_1),(0,v_2),\ldots,(0,v_n))\in (R \oplus I)^n$,
and  $\widetilde{w}=((0,w_1),(0,w_2),\ldots,(0,w_n)) \in (0 \oplus
I)^n$.  Hence as in the proof of Lemma \ref{key5}, we can write
$${\I}_n+XM(\widetilde{v},\widetilde{w})= \begin{cases} \underset{j=2}{\overset{n}\Pi}
\,\widetilde{\eps} \,ge_{1j}((0,Xw_{1j}))\,\widetilde{\eps}^{-1} &
\mbox{for linear case } \\ \underset{j\ne
2}{\underset{j=1}{\overset{n}\Pi}} \,\widetilde{\eps}
\,ge_{1j}((0,Xw_{1j}))\,\widetilde{\eps}^{-1} & \mbox{ otherwise }.
\end{cases} ~~~~~~~~~(\star)
$$

Now we split the proof into following two cases. \vp 

Case I: $\eps$ is an elementary generator of the type $ge_{pq}(x)$,
$x\in R$. First applying the homomorphism  $X\mapsto X^2$ and then
applying Lemma \ref{key3} over $R[X]$ we get $${\I}_n +
X^2M(\widetilde{v},\widetilde{w})=\underset{j}\Pi\left(\underset{t=1}{\overset{k}\Pi}
ge_{p_{j(t)}q_{j(t)}}(X\widetilde{h}_{j(t)}(X))\right),$$ where
$\widetilde{h}_{j(t)}(X)\in ((0\oplus I)[X])$.  Again, as before
applying the homomorphism $\phi$ it follows that
 $${\I}_n+X^2M(v,w)=\underset{j}\Pi\left(\underset{t=1}{\overset{k}\Pi}
ge_{p_{j(t)}q_{j(t)}}(Xh_{j(t)}(X))\right),$$ where $h_{j(t)}(X)\in
I[X]$; as desired.  Hence the result also follows for $d\gg 0$. \vp 

Case II: $\eps$ is a product of elementary generators of the type
$ge_{pq}(x)$. Let $\mu(\eps)=r$. Arguing as before, the result follows
by  applying the homomorphism $X\mapsto X^{2^r}$  using the Corollary
\ref{key4}.  \vp 

{\noindent}(5) $\Rightarrow$ (4): Given that $\alpha_s(X)\in
{\E}(n,R_s[X], I_s [X])$, where $s$ is non-nilpotent element in the
ring $R$,  and $\alpha(0)={\rm {\I}_{n}}$.  By Lemma \ref{elementary}, there 
exists $\widetilde{\alpha}_{(s,0)}(X)\in {\E}(n,R_s[X]\oplus I_s
[X])$,  where the element $(s,0)$ will remain non-nilpotent in the
ring $R\oplus I$, and $\phi(\widetilde{\alpha}_{(s,0)}(X))=\alpha_s(X)$. 

Also, by Lemma \ref{key1},
$\alpha_s(X)$ can be written as a product of the matrices of the form
$\eps_s ge_{ij}(Xh(X))\eps_s^{-1}$, with $h(X) \in I_s[X]$ and $\eps_s\in {\E}(n, R_s)$. 
Hence using the proof of Lemma \ref{elementary}  it
follows that $\widetilde{\alpha}_{(s,0)}(X)$ can be written as a
product of the matrices of the form
$\widetilde{\eps}_{(s,0)}ge_{ji}((0,Xh(X)))\widetilde{\eps}_{(s,0)}^{-1}$, where
$\phi(\widetilde{\eps}_{(s,0)})=\eps_s$ and $(0,Xh(X)) \in ((R\oplus
I)_{(s,0)}[X])$.

Applying the homomorphism $X\mapsto XT^d$, where $d\gg 0$, from the
polynomial ring $R[X]$ to the polynomial ring $R[X,T]$, we consider
$\widetilde{\alpha}_{(s,0)}(XT^d)$. Note that $R_s[X,T]\cong
(R_s[X])[T]$. Now, using the Equation ($\star$) as in the proof of
$(6)\Ra (5)$,  we can rewrite $\widetilde{\alpha}_{(s,0)}(XT^d)$ as
the form ${\I}_{n} + XT^d M(v,w)$; for some  suitable $v, w$ over the
ring $(R_s[X]\oplus I_s[X])[T]$. Hence by $(5)$  we can write
$\widetilde{\alpha}_{(s,0)}(XT^d)$ as a product of elementary
generators of general linear (symplectic/orthogonal resp.) group  such
that each of those elementary generator is congruent to  identity
modulo the ideal $(T)$  over the ring  $((R_s\oplus I_s)[X])[T]$.
Let $l$ be the maximum of the powers occurring in the denominators
of those elementary generators. Again, as $R$ assumed to be
Noetherian,  by applying the homomorphism $T\mapsto (s,0)^mT$, for
$m\ge l$, it follows from Lemma \ref{noeth} that by (uniquely)
identifying it's lift  over the ring  $(R\oplus I)[X,T]$  we can write
$\widetilde{\alpha}_{(s,0)}((s,0)^mXT^d)$ as a product of elementary
generators of the general linear (symplectic/orthogonal resp.) group
such that each of those elementary generator is congruent to identity
modulo $(T)$. {\it i.e.}  there exists some $\widetilde{\beta}(X,T)\in
{\E}(n,(R\oplus I)[X,T])$ such that $\widetilde{\beta}(0,0)={\rm
I}_{n}$ and
$\widetilde{\beta}_{(s,0)}(X,T)=\widetilde{\alpha}_{(s,0)}((b,0)XT^d)$
for some $(b,0)\in (s,0)^m(R\oplus I)$.  Finally, by substituting
$T=(1,0)$ and using Lemma \ref{noeth},  we get
$\widetilde{\alpha}((b,0)X)\in {\E}(n,(R\oplus I)[X])$. Hence the
result follows applying $\phi$ as before. \vp 

{\noindent}(4) $\Rightarrow$ (3): Since $\alpha_{\m}(X)\in
{\E}(n,R_{\m}[X], I_{\mf m}[X])$, for all $\m \in {\M}(R)$, for each
$\m$ there exists $s\in R \setminus \m$ such that  $\alpha_s(X)\in
{\E}(n,R_s[X], I_s[X]).$  Observe that $$R_s[X]\oplus I_s[X] \cong
(R_s\oplus I_s)[X]\cong (R\oplus I)_s[X].$$ Hence by Lemma
\ref{elementary}, applied to the base ring $R_s[X]$, there exists
$\widetilde{\alpha}_{(s,0)}(X)\in  {\E}(n,(R\oplus I)_{(s,0)}[X])$
such that $\phi_s(\widetilde{\alpha}_{(s,0)})=\alpha_s$.
Let $$\widetilde{\theta}(X,T)=
\widetilde{\alpha}_{(s,0)}(X+T)\widetilde{\alpha}_{(s,0)}(T)^{-1}.$$
Then $\widetilde{\theta}(X,T)\in {\E}(n,(R\oplus I)_{(s,0)}[X,T])$ and
$\widetilde{\theta}(0,T)={\I}_n$.  By the condition (4) of the Theorem,
applied to the base ring $(R\oplus I)[T]$, there exists
$\widetilde{\beta}(X)\in  {\E}(n,(R\oplus I)[X,T])$ such that 
 \begin{equation} \label{A} \widetilde{\beta}_{(s,0)}(X)=
\widetilde{\theta}((b,0)X,T).
 \end{equation} with $(b,0)\in (s,0)^l(R\oplus I)$ for some $l\gg 0$. 

Now, using the Noetherian property of $R\oplus I$, as mentioned in the
Remark \ref{finite}, we may consider  a finite cover of $R\oplus I$,
say $(s_1,0)+\cdots+ (s_r,0)=(1,0)$.   Since for $l\gg 0$, the ideal
$\langle (s_1,0)^l,\ldots, (s_r,0)^l\rangle = R\oplus I$, we choose
$(b_1,0),\ldots,(b_r,0)\in R\oplus I$, with $(b_i,0)\in
(s_i,0)^l(R\oplus I),\, l\gg 0$ such that (\ref{A}) holds and
$(b_1,0)+\cdots+(b_r,0)=(1,0)$.  Hence for each $i=1,\ldots,r$,  there
exists ${\widetilde{\beta}}^i(X)\in  {\E}(n,(R\oplus I)[X,T])$ such
that  $\widetilde{\beta}^i_{(s_i,0)}(X)=
\widetilde{\theta}((b_i,0)X,T)$.
Now, $$\underset{i=1}{\overset{r}\Pi}\widetilde{\beta}^i(X)\in
{\E}(n,(R\oplus I)[X,T]).$$ But, $$\widetilde{\alpha}_{s_1'\cdots
s_r'}(X)= \left(\underset{i=1}{\overset{r-1}\Pi}
\widetilde{\theta}_{s_1'\cdots \widehat{s_i'}\cdots s_r'}(b_i'X,T)
{\mid}_{T=b_{i+1}'X+\cdots +b_r'X}\right)
\widetilde{\theta}_{s_1'\cdots \cdots s_{r-1}'}(b_r'X,0),$$ where
$s_i'=(s_i,0)$ and $b_i'=(b_i,0)$ for each $i=1,\ldots,r$.  Now
$\widetilde{\alpha}(0)={\I}_n$. Also, as a consequence of the Lemma
\ref{noeth} it follows that  the map  $${\E}(n,R,(s,0)^k(R\oplus I)[X]) \ra
{\E}(n,(R\oplus I)_{(s,0)}[X])$$ for $k\in \mbb{N}$ is injective for each
$s=s_i$.  Hence by (uniquely) identifying $\widetilde{\alpha}_{s_1'\cdots
s_r'}(X)$ with its lift, we conclude  $\widetilde{\alpha}(X)\in {\E}(n,R[X]\oplus
I[X])$. Finally, applying the map $\phi$ we get  $\alpha(X)\in {\E}(n,R[X],
I[X])$; as desired. \vp  

{\noindent}(3) $\Rightarrow$ (2): This is the implication where we use
the commutative property of the base ring $R$.  Let $\alpha(X)={\I}_{n} +
XM(v,w)$, where  $v \in {\Um}_n(R,I)$ and $\langle v,w\rangle=0$ and $w
\in I^n$.  Then $\alpha(0)={\I}_n$. Let $v = (1+v_1,v_2, \ldots ,v_n)$
and $w=(w_1,w_2,\ldots,w_n)\in I^n$, with $v_i, w_i \in I$  for
$i=1,\ldots, n$. Then by Lemma \ref{key5}, $\alpha_{\mf{m}}(X)$ is elementary 
for every maximal ideal $\mf{m}$ in $R$. Hence $\alpha(X)$ is elementary by $(3)$. 
\vp\\

{\noindent}(2) $\Rightarrow$ (1): Let $\epsilon \in {\E}(n,R,I)$ and
$\gamma =(\gamma_{ij})\in {\G}(n,R)$. There exist $\widetilde{\epsilon} \in {\E}(n,R
\oplus I)$ and $\widetilde{\gamma} = ((\gamma_{ij},0)) \in {\G}(n,R \oplus I)$ respectively
such that $\phi (\widetilde{\epsilon} )=\epsilon$ and $\phi
(\widetilde{\gamma} ) =\gamma$.    Using (2) $\Rightarrow$  (1) of the
absolute case we get $\widetilde{\gamma}\,
\widetilde{\epsilon}\,\widetilde{\gamma}^{-1}\in   {\E}(n,R \oplus I)$ as
${\E}(n,R \oplus I) \triangleleft {\G}(n,R \oplus I)$ and applying the
homomorphism $\phi$   it follows ${\gamma}{\epsilon}{\gamma}^{-1}\in
{\E}(n,R,I)$; as required.\vp

{\noindent}(1) $\Rightarrow$ (7): Let $v=\gamma e_1$ where $\gamma \in
{\G}(n,R)$. Then there exists $\widetilde{\gamma} \in {\G}(n,R \oplus I)$
such that  $\phi (\widetilde{\gamma} ) =\gamma$. Let $\widetilde{v} =
\widetilde{\gamma} e_1$ and  $\widetilde{w} = ((0,w_1),(0,w_2),\ldots
,(0,w_n)) \in (0 \oplus I)^n$. We have $\langle
\widetilde{v},\widetilde{w}\rangle= 0$. Hence, using (1) $\Rightarrow$
(7) of the absolute case it follows that ${\I}_n
+M(\widetilde{v},\widetilde{w}) \in E_n (R \oplus I )$. Now  applying
the homomorphism $\phi$  we get ${\I}_n +M(v,w)\in {\E}(n,R,I)$; as
required. \vp 

The above implications prove the equivalence of the statements. \hb \vp

\begin{re} \tn{Assuming the result for the absolute case treated in \cite{rrr} one can give
simpler proofs of the steps (5) $\Rightarrow$ (4), and (4)
$\Rightarrow$ (3).  But, there is a gap in the proof of the absolute
case in \cite{rrr}, as mentioned in \cite{Basu}. The gap was filled 
in  \cite{Basu} by proving results for Bak's unitary groups, which
cover linear, symplectic and orthogonal groups, and some more
classical type groups.  To make this note self-contained,
we have given the detailed proofs of those steps. }
\end{re}

\section{Relative L-G Principle for Transvection Subgroups}

In this section we shall state auxiliary results without detailed
proofs. For definitions of  symplectic and orthogonal modules and their
transvection subgroups we refer to \cite{BBR}. 

{\noindent} In \cite{BBR}, the first and third authors together with Anthony
Bak generalized Quillen-Suslin's  local-global principle for the
transvection subgroups of the projective, symplectic and orthogonal
modules.  As before, all three cases were treated uniformly. 
We observe below how to obtain relative versions of that local-global principle.
To state the results we need to recall a few notations. 

\begin{nt} \label{note} \tn{In the sequel $P$ will denote either a
finitely generated  projective $R$-module of rank  $n$, a symplectic
$R$-module or an orthogonal $R$-module of even rank  $n=2m$ with a
fixed form $\langle \, ,\rangle$. And  $Q$ will denote $P\oplus R$ in
the linear case, and $P\perp R^2$, otherwise. 
We will use the notation $Q[X]$
to denote $(P\oplus
R)[X]$ in the linear case and  $(P\perp R^2)[X]$, otherwise. We assume that the rank of the projective module is 
$n\ge 2$,  when dealing with the linear case, and $n\ge 6$, when
considering the  symplectic and the orthogonal cases. For a finitely
generated projective  $R$-module $M$ we use the notation }
\tn{${\G}(M)$ to denote ${\aut}(M)$, ${\Sp}(M,\langle \, ,\rangle )$
and  ${\O}(M,\langle \, ,\rangle )$ respectively;  
denote ${\SL}(M)$, ${\Sp}(M,\langle \, ,\rangle )$ and
${\tran}(M)$, ${\tran}_{\Sp}(M)$ and  ${\tran}_{\O}(M)$ respectively;
and  $\tn{ET}(M)$ to denote ${\fl}(M)$, ${\fl}_{\Sp}(M)$ and
${\fl}_{\O}(M)$ respectively.}

\end{nt}  

We shall also assume the following hypotheses:  \vp\\ (H1) for every maximal
ideal $\mf{m}$ of $R$, the symplectic (orthogonal)  module
$Q_{\mf{m}}$ is isomorphic  to $R^{2m+2}_{\mf{m}}$ for the standard
bilinear form $\mbb{H}(R_{\mf{m}}^{m+1})$.  \vp\\ (H2) for every
non-nilpotent $s\in R$, if the projective module  $Q_s$ is free
$R_s$-module, then the symplectic (orthogonal) module  $Q_s$ is
isomorphic to $R^{2m+2}_s$ for the standard bilinear form
$\mbb{H}(R_s^{m+1})$.

We recall the following fact just to remind the reader that in the
free case the  transvection subgroups coincides with the elementary
subgroups. Here the maps  $\varphi$, $\varphi_p$, $\sigma$ and $\tau$
are as defined in \cite{BBR}. 
\begin{lm} \label{free1} If the projective module $P$ is free of
finite  rank $n$\, $($in the symplectic and the orthogonal
cases we assume that the projective module is free for the standard
bilinear  form$)$, then  ${\tran}(P)={\E}_{n}(R)$,  ${\tran}_{\Sp}(P)=
{\ESp}_{n}(R)$ and  ${\tran}_{\O}(P)= {\EO}_{n}(R)$ for $n\ge 3$ in
the  linear case and for $n\ge 6$ otherwise. 
\end{lm}  {\bf Proof.}  In the linear case, for $p\in P$ and
$\varphi\in P^{*}$ if  $P=R^n$ then $\varphi_{p} : R^{n}\ra R\ra
R^{n}$. Hence  $1+\varphi_{p}={\I}_{n}+v.w^t$ for some column vectors $v$ and
$w$ in $R^{n}$. Since $\varphi(p)=0$, it follows that
$\langle v,w\rangle=0$. Since either $v$ or $w$ is unimodular,  it
follows that 
$1+\varphi_{p}={\I}_{n}+v.w^t \in{\E}_{n}(R)$.  Similarly, in the
non-linear cases we have
$\sigma_{(u,v)}(p)={\I}_{n}+v.\widetilde{w}+w.\widetilde{v}$, and
$\tau_{(u,v)}(p)={\I}_{n}+v.\widetilde{w}-w.\widetilde{v}$, where  either
$v$ or $w$ is unimodular and $\langle v,w\rangle=0$. (Here
$\sigma_{(u,v)}$ and $\tau_{(u,v)}$ are as in the definition of
symplectic  and orthogonal transvections.) Classically, these are
known to be elementary  matrices - for details see \cite{SUS} for the
linear case, \cite{KOP} for the  symplectic case, and \cite{SUSK} for
the orthogonal case. \hb 
\hb
\begin{re} \tn{Lemma \ref{free1} holds for $n=4$ in  the symplectic
case. This will follow from Remark  \ref{free2}.}
\end{re}
\begin{re} \label{free2}  \tn{${\ESp}_4(A)$ is a normal subgroup of
${\Sp}_4(A)$ by  (\cite{KOP}, Corollary 1.11). Also  ${\ESp}_4(A[X])$
satisfies the Dilation Principle and the Local-Global  Principle by
(\cite{KOP}, Theorem 3.6).  Since we were intent on a uniform proof,
these cases  have not been covered above by us.}
\end{re} 

\begin{pr} \label{dil} {\bf (Relative Dilation Principle)}  Let $R$ be
a commutative ring with identity, and $I\subsetneq R$ an ideal in $R$.
Let $P$ and $Q$ be as in \ref{note}. Assume that $($H2$)$ holds.  Let
$s$ be a non-nilpotent in $R$ such that $P_s$ is free, and  let
$\sigma(X)\in {\G}(Q[X], I[X])$ with $\sigma(0)=\tn{Id}$. Suppose 
$$\sigma_s(X)\in \begin{cases} {\E}(n+1,R_s[X], I_s[X]) & \mbox{in the linear case, }\\
{\E}(n+2,R_s[X], I_s[X])& \mbox{ otherwise. } 
\end{cases}$$  Then there exists $\widehat{\sigma}(X)\in {\rm
ET}(Q[X], I[X])$  and $l>0$ such that $\widehat{\sigma}(X)$ localizes
to  $\sigma(bX)$ for some $b\in (s^l)$ and
$\widehat{\sigma}(0)=\tn{Id}$.
\end{pr}  {\bf Proof.} Follows by imitating the technique explained in
\cite{BBR}, and following steps mentioned in Section 2. \hb 
\begin{tr} {\bf (Relative Local-Global Principle)} \label{lgt} Let $R$
be a commutative ring with identity, and $I\subsetneq R$ an ideal in
$R$.  Let $P$ and $Q$ be as in \ref{note}. Assume that $($H1$)$ holds.
Suppose $\sigma(X)\in {\G}(Q[X], I[X])$ with $\sigma(0)=\tn{Id}$. If 
$$\sigma_{\mf{p}}(X)\in \begin{cases} {\E}(n+1,R_{\mf{p}}[X],I_{\mf{p}}[X]) & 
\mbox{in the linear case, }\\ {\E}(n+2,R_{\mf{p}}[X],I_{\mf{p}}[X])&
\mbox{ otherwise } 
\end{cases}$$  for all $\mf{p}\in {\sp}(R)$, then $\sigma(X)\in {\rm
ET}(Q[X],I[X])$.
\end{tr}  {\bf Proof.} Follows by using similar technique as in
$(4)\Ra (3)$ in Theorem 3.1 of \cite{rrr}, and  and arguing as in
Section 2. \hb

\begin{re} \tn{The authors believe that  the above method using
the ``Noetherian Excision ring'', makes it possible to deduce 
the relative versions of almost all the results mentioned in \cite{rrr},
\cite{BBR}, \cite{Basu},  and the results mentioned in \cite{RB1}
between pgs 35--40. }
\end{re}

 \vspace{0.05cm}

{\small

\noindent
{\it Indian Institute of
Science Education and Research, Dr. Homi Bhabha Road, Pune 400 008, India}.\\ 
{\it email: rabeya.basu@gmail.com, rbasu@iiserpune.ac.in}\vp 

\noindent
{\it  Somaiya College, Vidyavihar, Mumbai 400 077, India.}\\ 
{\it email: reemag16@gmail.com} \vp

\noindent
{\it Tata Institute of Fundamental Research, 1, Dr. Homi Bhabha Road,
Mumbai \lb 400 005, India.}\\ 
{\it email: ravi@math.tifr.res.in} }

\end{document}